\documentstyle{amsart}

\numberwithin{equation}{section}

\newtheorem{theorem}{Theorem}[section]

\newtheorem{lemma}[theorem]{Lemma}

\newtheorem{proposition}[theorem]{Proposition}

\theoremstyle{remark}

\begin{document}

\title{Entire invariant solutions to Monge-Amp\`ere equations}

\author{Roger Bielawski}
\thanks{Research supported by an EPSRC advanced fellowship}

\address{Department of Mathematics, University of Glasgow, Glasgow G12 8QW, UK }

\email{R.Bielawski@@maths.gla.ac.uk}

\begin{abstract}We prove existence and regularity of entire solutions to
Monge-Amp\`ere equations invariant under an irreducible action of a compact
Lie group.
 \end{abstract}

\maketitle

\thispagestyle{empty}

 We consider Monge-Amp\`ere equations of the form
\begin{equation} f(\nabla
\phi)\det D_{ij}\phi=g(x)\label{MA} \end{equation}
 where $f$ and $g$ are nonnegative measurable functions on ${\Bbb R}^n$.
 We recall first
the concept of a weak solution of \eqref{MA}. Let $\phi$ be a convex
function. Then $\nabla \phi$ is a well-defined multi-valued mapping:
$(\nabla\phi)(x)$ is the set of slopes of all supporting hyperplanes to the
graph of $\phi$ at $(x,\phi(x))$. If $B$ is a subset of ${\Bbb R}^n$, let
$\nabla\phi(B)$ be its image in the multi-valued sense. Then $\phi$ is a
weak solution of \eqref{MA} if
\begin{equation} \int_B g(x)dx=\int_{\nabla\phi(B)}f(y)dy \label{weak} \end{equation}
for every Borel set $B$. Let us denote the right-hand side by
$\omega(B,\phi,f)$. It can be shown that it is a Borel measure on ${\Bbb
R}^n$. A basic result \cite{Bak} is that if $u_k\rightarrow u$ compactly and
$f_k \rightarrow f$ uniformly, then $\omega(\cdot,u_k,f_k)$ converges to
$\omega(\cdot,u,f)$ weakly, i.e. as functionals on the space of compactly
supported continuous functions.
\par
 Existence of
weak solutions to \eqref{MA} defined on all of ${\Bbb R}^n$ have been shown
in \cite{Bak} under the assumption that $g\in L^1({\Bbb R}^n)$ and in
\cite{CW} in the case when $f=1$ and $g,1/g$ are bounded. In this note we
wish to give a simple proof of existence and regularity of entire solutions
to \eqref{MA} under a different type of assumption: $f$ and $g$ are
invariant under an irreducible action of a compact (Lie) group. More
precisely we shall prove:
\begin{theorem} Let $K\subset O(n)$ be a compact subgroup acting
irreducibly on ${\Bbb R}^n$ and let $f,g$ be two  nonnegative $K$-invariant
 measurable functions on ${\Bbb R}^n$. Furthermore, assume that $f$ and $g$ are locally
bounded and that $$\int_{{\Bbb R}^n}f=+\infty.$$ Then there exists a (weak)
convex $K$-invariant solution $\phi:{\Bbb R}^n\rightarrow {\Bbb R}$ of the
Monge-Amp\`ere equation \eqref{MA}.\label{Mat}
\end{theorem}

\begin{theorem}
In the situation of Theorem \ref{Mat} suppose, in addition,  that $f$ and
$g$ are strictly positive and of class $C^{p,\alpha}$, $p\geq 0$, $\alpha>0$. Then any entire
$K$-invariant convex weak solution $\phi$ of \eqref{MA} is
$C^{p+2,\alpha}$.\label{smooth}
\end{theorem}
Well-known examples show that this regularity result does not hold without
the assumption of $K$-invariance, even when $f=g=1$.
\par
The above theorems have been proved before in \cite{rb} in the case when $K$
is a reflection group. They can be applied to show existence of K\"ahler
metrics with prescribed Ricci curvature on complexified symmetric spaces.

\section{Proofs}

We begin with a simple lemma:
\begin{lemma}
Let $K\subset O(n)$ be a compact subgroup acting irreducibly on ${\Bbb
R}^n$. For any $x\in {\Bbb R}^n$, let $C_x$ be the convex hull of the orbit
$K\cdot x$. There exists an $\epsilon >0$ such that the ball $B(0,\epsilon)$
is contained in $C_x$ for any $x$ with $|x|=1$. \label{irrep}
\end{lemma}
\begin{pf}
Suppose that the conclusion does not hold. Then there exists a sequence of
points $x_k$ with $|x_k|=1$ such that the boundary of $C_{x_k}$ contains a
point $y_k$ with  $y_k\rightarrow 0$. By choosing a convergent subsequence
of $x_k$ we find a point a point $x_\infty$ on the unit sphere such that $0$
lies on the boundary of $C_{x_\infty}$. As $C_{x_\infty}$ is convex, it lies
to one side of any of its supporting hyperplanes at $0$. Hence the orbit
$K\cdot x_\infty$ lies to one side of a hyperplane $H$ passing through $0$.
Therefore the center of mass of the orbit $$ m=\frac{1}{|K|}\int_K g\cdot
x_\infty dg$$ cannot be the origin unless the whole orbit lies on $H$. On
the other hand $m$ has to be the origin as it is a $K$-invariant point and
the representation is irreducible. It follows that the orbit $K\cdot
x_\infty$ is contained in a hyperplane and from this the subspace generated
linearly by $K\cdot x_\infty$ is a proper $K$-invariant subspace of ${\Bbb
R}^n$, a contradiction.
\end{pf}

We now prove Theorem \ref{Mat}. Put $f_k=f+1/k$ and $g_k=g+1/k$, $k\in {\Bbb
N}_+$. Let $B_r$ denote the ball of radius $r$ centred at the origin. Let
$R_k$ be a number defined by $$\int_{B_{R_k}}f_k=\int_{B_k}g_k.$$ According
to \cite{Bre,Cafconvex} there exists a unique (up to a constant) strictly
convex solution $\phi_k$ of
\begin{equation} f_k(\nabla \phi_k)\det D_{ij}\phi_k=g_k(x)\label{MAk}
\end{equation}
which is of class $C^{1,\beta}$, for some $\beta>0$, and such that $\nabla
\phi_k$ maps $B_k$ onto $B_{R_k}$. Moreover, as $K\subset O(n)$ and all the
data are $K$-invariant, $\phi_k$ is $K$-invariant (this follows from
uniqueness, since $\phi_k\circ g$, $g\in K$, is also a solution of
\eqref{MAk}). To prove the existence of a weak solution to \eqref{MA},
defined on all of ${\Bbb R}^n$, it is enough to show that the functions
$\phi_k$ are uniformly (in $k$) bounded on any ball $B_R$, $k\geq R$.
Indeed,
 a bounded sequence of convex functions on a bounded open convex domain has
a convergent subsequence. This follows from the elementary fact, which we
will use repeatedly, that the slopes of supporting hyperplanes of a convex
function, bounded by $R$ on a domain $G$, are bounded by $2R/\delta$ on any
subdomain $G^\prime$ such that $\text{dist}(G^\prime,\partial G)\geq \delta$.
\par
Let us show that the functions $\phi_k$, $k\geq R$, $\phi_k(0)=0$, are
bounded on $B_R$ uniformly in $k$. It is enough to show that $\nabla \phi_k$
are bounded uniformly on  $B_R$. Suppose that there exists a sequence of
points $x_{k_j}$ in $B_R$ such that $\bigl|\nabla
\phi_{k_j}(x_{k_j})\bigr|\geq j$. Let $y_{k_j}=\nabla \phi_{k_j}(x_{k_j})$.
Consider the unique solution $\psi_k$, $\psi_k(0)=0$, to the equation $$
g_k(\nabla \psi_k)\det D_{ij}\psi_k=f_k(y)$$ mapping $B_{R_{k_j}}$ onto
$B_{k_j}$. According to \cite{Bre}, $\nabla \psi_k$ is the inverse of
$\nabla \phi_k$. Thus $\nabla \psi_{k_j}(y_{k_j})=x_{k_j}$ and hence
$\bigl|\nabla \psi_{k_j}(y_{k_j})\bigr|\leq R$. Observe that $\nabla
\psi_k(0)=0$ since $\nabla \psi_k$ is $K$-equivariant and the representation
is irreducible. Let $v_{k_j}$ be the unit vector in the direction $y_{k_j}$.
Since $\psi_k$ is convex, $\langle\nabla\psi_{k_j}(y), v_{k_j}\rangle\leq R$
for all $y\in [0,y_{k_j}]$ and therefore $\bigl|\psi_{k_j}(y)\bigr|\leq
R|y|$ for such $y$. Using now Lemma \ref{irrep} we conclude that for any
$r>0$, $\psi_{k_j}$ is bounded by $\frac{R}{\epsilon}r$ on $B_r$ for large
$j$. It follows that there is a subsequence of $\psi_{k_j}$ convergent to a
convex solution $\psi:{\Bbb R}^n\rightarrow {\Bbb R}$ of $$g(\nabla\psi)\det
D_{ij}\psi=f.$$ The function $\psi$ is bounded by $Rr/\epsilon$ on any
$B_r$. Since $\psi$ is convex, $\nabla\psi$ is bounded by $4R/\epsilon$ on
$B_{r/2}$. Thus $U=\nabla\psi\bigl({\Bbb R}^n\bigr)$ is contained in
$B_{4R/\epsilon}$. However $$+\infty=\int_{{\Bbb R}^n}f= \int_U g, $$ which
leads to a contradiction. Theorem \ref{Mat} is proved.

\begin{proposition}  Suppose, in addition, that $\int_{{\Bbb R}^n}g>0$. Then
any entire $K$-invariant convex solution $\phi$ of \eqref{MA}  is a proper
function. \label{proper}
\end{proposition}
\begin{pf}
Suppose that there exists an unbounded sequence of points $x_k$ such that
$\phi(x_k)$ is bounded. Using convexity of $\phi$ and Lemma \ref{irrep} we
conclude that $\phi$ is a bounded function. This however implies that
$\phi$, being convex and defined on all of ${\Bbb R}^n$, is constant, which
contradicts the assumption on $g$.
\end{pf}

This proposition allows us to prove the regularity results of Theorem
\ref{smooth}. The proof is essentially the same as the one given by
Caffarelli and Viaclovsky \cite{CV}; we give it here for completeness.
\par
Consider the convex (and bounded by the last proposition)
sets $\Omega_c=\{x;\phi(x)\leq c\}$, $c\in {\Bbb R}$. For every $c$, we can find $R(c)>c$ such that $\text{dist}\,(\Omega_c,\partial\Omega_{R(c)})\geq 1$. By the argument used before, the slopes $\nabla\phi$ of $\phi$ are bounded by $2R(c)$ on $\Omega_c$. If $f$ and $g$ are strictly positive and continuous,  then equation \eqref{weak} implies  that there are constants $\lambda_1,\lambda_2$, depending on $c$, such that
\begin{equation} \lambda_1|B|\leq |\nabla\phi(B)|\leq \lambda_2|B| \label{inequality}\end{equation}
for any Borel subset $B\subset \Omega_c$.
We can apply Corollary 2 in
\cite{Caflocalization} to $\Omega_{c+1}$ and conclude that $\phi$ is strictly convex at any point of $\Omega_c$. Since $c$ was arbitrary, $\phi$ is strictly convex everywhere. Now the main result of \cite{Cafsome} implies
that $\phi$ is locally of class $C^{1,\beta}$, where ``locally" means that $\beta$ can be chosen uniformly only on bounded subsets of ${\Bbb R}^n$.
We observe that the only additional properties of $f$ and $g$ we have used so far is the boundedness of $1/f$ and $1/g$. Therefore we have proved:
\begin{proposition} In the situation of Theorem \ref{Mat} suppose, in addition, that $1/f$ and $1/g$ are locally bounded functions and let $\phi$ be  an entire $K$-invariant convex weak solution of \eqref{MA}. Then $\phi$ is of class $C^{1,\beta}$ on any compact subset $P$, where $\beta$ may depend on $P$.\hfill $\Box$\end{proposition}
To finish the proof of  Theorem \ref{smooth}, assume that $f$ and $g$ are of class $C^{0,\alpha}$. Now,  $\phi$ is a solution of $\det D_{ij}\phi=\tilde g$, where
$\tilde{g}=g/f(\nabla \phi)$. Since, by the last proposition,  $\nabla \phi$ is H\"older continuous,  $\tilde{g}$ is  (locally) of class $C^{0,\gamma}$. Using \eqref{inequality}, we can  apply
Theorem 2 in \cite{Cafinterior} to $u(x)=\phi(x) -c$, $c\in {\Bbb R}$ and
the convex and bounded set $\Omega_c=\{x; u(x)\leq 0\}$ to conclude that
$\phi$ is locally of class $C^{2,\gamma}$ and hence globally $C^{2}$. Therefore $\tilde{g}=g/f(\nabla \phi)$ is of class   $C^{0,\alpha}$ everywhere and, repeating the argument, we conclude that $\phi$ is $C^{2,\alpha}$.
Higher regularity is standard.


\begin{thebibliography}{99}

\bibitem{Bak}
{I.J. Bakelman} {\sl Convex analysis and nonlinear geometric elliptic
equations}, Springer-Verlag, Berlin, 1994.


\bibitem{rb}
{R.Bielawski} `K\"ahler metrics on $G^{\Bbb C}$', {\it J. Reine Angew.
Math.}, to appear.

\bibitem{Bre}
{Y. Brenier} `Polar factorization and monotone rearrangement of vector-valued
functions', {\it Comm. Pure Appl. Math.} 44 (1991),  375--417.


\bibitem{Caflocalization}
{L.A. Caffarelli} `A localization property of viscosity solutions to the
Monge-Amp\`ere equation and their strict convexity', {\it Ann. of Math. (2)}
131 (1990), 129--134.


\bibitem{Cafinterior}
{L.A. Caffarelli} `Interior $W\sp {2,p}$ estimates for solutions of the
Monge-Amp\`ere equation', {\it Ann. of Math. (2)} 131 (1990), 135--150.


\bibitem{Cafsome}
{L.A. Caffarelli} `Some regularity properties of solutions of Monge-Amp\`ere
equation', {\it Comm. Pure Appl. Math.} 44 (1991), 965--969.


\bibitem{Cafconvex}
{L.A. Caffarelli} `The regularity of mappings with a convex potential', {\it
J. Amer. Math. Soc.} 5 (1992), 99--104.

\bibitem{CV}
{L.A. Caffarelli \and J.A. Viaclovsky} `On the regularity of solutions to
Monge-Ampère equations on Hessian
   manifolds', {\it Comm. Partial Differential Equations} 26 (2001), no. 11-12, 2339--2351.


\bibitem{CW}
{K.-S. Chou \and X.-J. Wang} `Entire solutions of the Monge-Amp\`ere
equation', {\it Comm. Pure Appl. Math.} 49 (1996),  529--539.






\end{thebibliography}
\end{document}